\documentclass[a4paper,12pt]{amsart}
\usepackage{amsmath}
\usepackage{amssymb}
\usepackage{amsthm}
\usepackage{mathrsfs}
\numberwithin{equation}{section}
\theoremstyle{plain}
\newtheorem{thm}{Theorem}[section]
\newtheorem{cor}[thm]{Corollary}
\newtheorem{lem}[thm]{Lemma}
\newtheorem{prop}[thm]{Proposition}
\theoremstyle{remark}

\begin{document}
\title[A Fourier restriction theorem]{A Fourier restriction theorem for hypersurfaces which are graphs of certain real polynomials}
\author[K.~Morii]{Kei MORII}
\address{Mathematical Institute, Tohoku University, Sendai, 980-8578, Japan.}
\email{sa3m28@math.tohoku.ac.jp}
\thanks{The author is supported by JSPS Research Fellowship for Young Scientists}
\begin{abstract}
We will extend the Fourier restriction inequality 
for quadratic hypersurfaces obtained by Strichartz. 
We will consider the case where the hypersurface 
is a graph of a certain real polynomial 
which is a sum of one-dimensional monomials. 
It is essential to examine the decay of a one-dimensional oscillatory integral.
\end{abstract}
\maketitle
%
%
\section{Introduction}
Let $S$ be a hypersurface in $\mathbb{R}^n$, 
$n\geqslant 2$. 
We consider the Fourier restriction inequality 
\begin{equation}
\left(
\int_S \lvert\hat \phi(\xi)\rvert^2 d\mu_n (\xi)
\right)^{1/2}
\leqslant 
C_p 
\lVert\phi\rVert_{L^p(\mathbb{R}^n)}
\text{\,\,for\,\,}
\phi \in L^p(\mathbb{R}^n), 
\label{restriction_ineq}
\end{equation}
where the measure $d\mu_n$ on $S$ is defined as follows: 
$$
d\mu_n(\xi)
=\left\lvert
\frac{\partial\tilde R}{\partial \xi_n}(\xi)
\right\rvert^{-1}
d\xi_1\ldots d\xi_{n-1}
$$
when $S$ is written as 
$S=\{\xi\in\mathbb{R}^n;\,\tilde R(\xi)=r\}$ 
with a constant $r\in \mathbb{R}$ 
and a real-valued function 
$\tilde R\in C^0(\mathbb{R}^n)$ 
which is partially differentiable 
with respect to $\xi_n$ and 
$(\partial\tilde R/\partial \xi_n)(\xi)\neq 0$ 
for almost every $\xi\in\mathbb{R}^n$. 
For $\xi=(\xi_1,\ldots,\xi_n)\in\mathbb{R}^n$, 
we write $\xi'=(\xi_1,\ldots,\xi_{n-1})$. 
In particular, 
if $S$ is the form of 
$S=\{\xi\in\mathbb{R}^n;\,\xi_n=R(\xi')\}$, 
then \eqref{restriction_ineq} becomes
\begin{equation}
\left(
\int_{\mathbb{R}^{n-1}} 
\lvert\hat \phi(\xi',R(\xi'))\rvert^2 d\xi'
\right)^{1/2}
\leqslant 
C_p 
\lVert\phi\rVert_{L^p(\mathbb{R}^n)}
\text{\,\,for\,\,}
\phi \in L^p(\mathbb{R}^n). 
\label{eq:restriction_ineq_rewritten}
\end{equation}
\par
For $p\in[1,\infty)$ and 
a subset of a Euclidean space $\Omega$, 
set
$$
\lVert f\rVert_{L^p(\Omega)}
=
\left(
\int_\Omega 
\lvert f(x)\rvert^p dx
\right)^{1/p}
$$
and let $L^p(\Omega)$ denote the set of 
all Lebesgue measurable functions $f$ on $\Omega$ 
such that $\lVert f\rVert_{L^p(\Omega)}<\infty$.
Let $i$ always denote the imaginary unit.
We define the Fourier transform 
in $x\in\mathbb{R}^n$ 
and the inverse Fourier transform 
in $\xi\in\mathbb{R}^n$ by setting
\begin{align*}
  \hat f(\xi)
& =
  (2\pi)^{-n/2}
  \int_{\mathbb{R}^n}
  f(x)e^{-ix\cdot\xi}dx,
\\
  \check{f}(x)
& =
  (2\pi)^{-n/2}
  \int_{\mathbb{R}^n}
  f(\xi)e^{ix\cdot\xi}d\xi, 
\end{align*}
respectively. 
Those of a generalized function 
are also denoted by the same notation. 
\par
In \cite{Strichartz1},
Strichartz determined 
the optimal range of the exponent $p$ 
for which \eqref{restriction_ineq} holds 
for all quadratic hypersurfaces $S$. 
A nondegenerate quadratic hypersurface $S$ 
which is not contained in an affine hyperplane 
is transformed into one of the following three types 
under an affine transformation. 
\begin{enumerate}
\item 
$S=\{\xi\in\mathbb{R}^n;\,
\xi_n=\xi_1^2+\cdots+\xi_s^2
-\xi_{s+1}^2-\cdots-\xi_{n-1}^2\}$, 
where $1\leqslant s\leqslant n-1$.
\item 
$S=\{\xi\in\mathbb{R}^n;\,
\xi_1^2+\cdots+\xi_s^2
-\xi_{s+1}^2-\cdots-\xi_n^2=0\}$, 
where $1\leqslant s\leqslant n-1$.
\item 
$S=\{\xi\in\mathbb{R}^n;\,
\xi_1^2+\cdots+\xi_s^2
-\xi_{s+1}^2-\cdots-\xi_n^2=1\}$, 
where $1\leqslant s\leqslant n$.
\end{enumerate}
The results obtained by Strichartz 
for the first type 
are the following.
\begin{thm}[Strichartz {\cite[Theorem~1]{Strichartz1}}]
\label{thm:tomoko}
Let $n\geqslant 2$, 
and let $S$ be a hypersurface 
$S=\{\xi\in\mathbb{R}^n;\,
\xi_n=\xi_1^2+\cdots+\xi_s^2
-\xi_{s+1}^2-\cdots-\xi_{n-1}^2\}$, 
where $1\leqslant s\leqslant n-1$.
Then, 
\eqref{eq:restriction_ineq_rewritten} holds with $C_p$ 
independent of $f$ if and only if
$$
p=\frac{2(n+1)}{n-3}. 
$$
\end{thm}
Strichartz also gives an estimate of solutions 
to the inhomogeneous Schr\"odinger evolution equations 
as an application of the Fourier restriction theorem. 
\begin{thm}[Strichartz {\cite[Corollary~1]{Strichartz1}}]
\label{thm:ayako}
Let $n\geqslant 1$. 
Let
$$
p
=\frac{2(n+2)}{n+4}, 
$$
and assume $\phi\in L^2(\mathbb{R}^n)$ 
and $f\in L^p(\mathbb{R}^{1+n})$. 
Let $u(t,x)$ be a solution to 
the initial value problem for
inhomogeneous partial differential equations
\begin{alignat*}{2}
D_t u+\Delta u
&=f(t,x) 
& &\text{\,\,in\,\,}\mathbb{R}^{1+n}, 
\\
u(0,x)
&=\phi(x) 
& &\text{\,\,in\,\,}\mathbb{R}^n, 
\end{alignat*}
where 
$$
D_t=-i\frac{\partial }{\partial t},
\,\,
\Delta
=
\sum_{j=1}^n
\frac{\partial^2}{\partial x_j^2}. 
$$
Then, 
$$
\lVert u\rVert_{L^{p/(p-1)}(\mathbb{R}^{1+n})}
\leqslant 
C(\lVert\phi\rVert_{L^2(\mathbb{R}^n)}
+\lVert f\rVert_{L^p(\mathbb{R}^{1+n})})
$$
holds with $C$ independent of $\phi$, 
$f$ and $u$. 
\end{thm}
The purpose of this paper is to study generalization of 
Theorems~\ref{thm:tomoko} and \ref{thm:ayako}. 
We will consider the case where $S$ 
is the graph of a certain real polynomial. 
We will introduce a method due to Strichartz, 
and extend Theorem~\ref{thm:tomoko} 
to prove the following. 
\begin{thm}
\label{thm:main_thm}
Let $n\geqslant 2$, 
and let $S$ be a hypersurface 
$S=\{\xi\in\mathbb{R}^n;\,\xi_n=R(\xi')\}$, 
where 
\begin{equation}
R(\xi')
=\sum_{j=1}^{n-1}a_j\xi_j^{k_j}, 
\label{eq:Rxi_prime}
\end{equation}
$a_j\in\mathbb{R}\setminus\{0\}$ and $k_j\in\{2,3,4,\ldots\}$ for all $j=1,\ldots,n-1$. 
Then, 
\eqref{eq:restriction_ineq_rewritten} holds with $C_p$ 
independent of $f$ if and only if
$$
p=
2-\frac{2}{2+\displaystyle\sum_{j=1}^{n-1}\frac{1}{k_j}}. 
$$
\end{thm}
The following theorem is our corresponding application 
of Theorem~\ref{thm:main_thm}.
\begin{thm}
Let $n\geqslant 1$. 
Let
$$
a(\xi)=\sum_{j=1}^n a_j\xi_j^{k_j}, 
$$
where $a_j\in\mathbb{R}\setminus\{0\}$ 
and $k_j\in\{2,3,4,\ldots\}$ for all $j=1,\ldots,n$. 
Let
$$
p
=2-\frac{2}{2+\displaystyle\sum_{j=1}^n \frac{1}{k_j}}, 
$$
and assume $\phi\in L^2(\mathbb{R}^n)$ 
and $f\in L^p(\mathbb{R}^{1+n})$. 
Let $u(t,x)$ be a solution to 
the initial value problem for
inhomogeneous partial differential equations
\begin{alignat}{2}
D_t u-a(D)u
&=f(t,x) 
& &\text{\,\,in\,\,}\mathbb{R}^{1+n}, 
\label{eq:pde_1}
\\
u(0,x)
&=\phi(x) 
& &\text{\,\,in\,\,}\mathbb{R}^n, 
\label{eq:pde_2}
\end{alignat}
where 
$$
D_t=-i\frac{\partial }{\partial t}
,\,\,
D=(D_1,\ldots,D_n),
\,\,
D_j=-i\frac{\partial }{\partial x_j}. 
$$
Then, 
$$
\lVert u\rVert_{L^{p/(p-1)}(\mathbb{R}^{1+n})}
\leqslant 
C(\lVert\phi\rVert_{L^2(\mathbb{R}^n)}
+\lVert f\rVert_{L^p(\mathbb{R}^{1+n})})
$$
holds with $C$ independent of $\phi$, 
$f$ and $u$. 
\label{thm:main_cor}
\end{thm}
\par
An essential matter in proving Theorem~\ref{thm:main_thm} 
is to examine the decay of the Fourier transform 
of $\exp(itR(x'))$, 
that is, 
\begin{equation}
\int_{\mathbb{R}^{n-1}}
\exp(ix'\cdot\xi'+itR(\xi'))d\xi'
\label{eq:os_int_n-1dim}
\end{equation}
for large $t$. 
Since $R$ is a sum of one-dimensional monomials, 
an estimate of \eqref{eq:os_int_n-1dim} 
is effectively reduced to that 
of a one-dimensional oscillatory integral
\begin{equation}
\int_{-\infty}^\infty
\exp(ix\xi+it\xi^k)d\xi,
\label{eq:os_int_1dim}
\end{equation}
where $k_j\in\{2,3,4,\ldots\}$. 
We will bound \eqref{eq:os_int_1dim} 
by $12\lvert t\rvert^{-1/k}$.
If the region of integration is a bounded interval, 
then
$$
\left\lvert
\int_\alpha^\beta \exp(ix\xi+it\xi^k)d\xi
\right\rvert
\leqslant 
C_{\alpha,\beta}\lvert t\rvert^{-1/k}
$$
immediately follows from the Van der Corput lemma. 
\par
The organization of this paper is as follows. 
In Section~\ref{sect:restriction}, 
we estimate the one-dimensional oscillatory integral \eqref{eq:os_int_1dim}, 
and give a proof of Theorem~\ref{thm:main_thm}. 
Section~\ref{sect:pde} 
describes a proof of Theorem~\ref{thm:main_cor}. 
\section{A Fourier restriction theorem}
\label{sect:restriction}
In this section, 
we estimate the one-dimensional oscillatory integral \eqref{eq:os_int_1dim}, 
and give a proof of Theorem~\ref{thm:main_thm}. 
We use a method due to Strichartz. 
Let $\Gamma$ denote the Gamma function. 
For $s\in\mathbb{R}$, 
let $s_+$ denote its positive part: 
$s_+=\max\{s,0\}$.
\begin{prop}[Strichartz {\cite[Lemma~2]{Strichartz1}}]
Let $n\geqslant 2$. 
We assume that a hypersurface $S$ is written as 
$S=\{\xi\in\mathbb{R}^n;\,\tilde R(\xi)=r\}$ 
with a constant $r\in \mathbb{R}$ and a real-valued function 
$\tilde R\in C^0(\mathbb{R}^n)$. 
Let
$$
G_z(\xi)
=
\frac{(\tilde R(\xi)-r)_+^z}{\Gamma(z+1)}
$$
for $z\in\mathbb{C}$. 
Moreover, 
we assume that for some $\lambda>1$, 
$\check G_{-\lambda+i\eta}$ is bounded: 
$$
\lVert\check G_{-\lambda+i\eta}\rVert_{L^\infty(\mathbb{R}^n)}
\leqslant C_\eta, 
$$
and that there exists $b<\pi$ such that
$$
\sup_{\eta\in\mathbb{R}}
e^{-b\lvert\eta\rvert}
\log C_\eta
<
\infty. 
$$
Then, 
\textup{\eqref{restriction_ineq}} holds for
$$
p
=\frac{2\lambda}{\lambda+1}. 
$$
\label{prop:Strichartz_method}
\end{prop}
\par
We use the following formula of integration later. 
\begin{lem}[{\cite[p. 360]{temp001}}]
$$
\frac{(2\pi)^{1/2}}{\Gamma(z+1)}
(\xi_+^z)^\wedge (x)
=
\lim_{\varepsilon\searrow 0}
\int_0^\infty 
\frac{e^{-\varepsilon \xi}\xi^z e^{ix\xi}}{\Gamma(z+1)}d\xi
=ie^{iz\pi/2}(x+i0)^{-z-1}
$$
for all $z\in\mathbb{C}$. 
\label{lem:formula_01}
\end{lem}
\par
Next, 
for $k_j\in\{2,3,4,\ldots\}$, 
we define a one-dimensional oscillatory integral
$$
A_k(x)
=
\lim_{\varepsilon\searrow 0}
\int_{-\infty}^\infty 
\exp(-\varepsilon\lvert s\rvert^k+ixs+is^k)ds
$$
for $x\in\mathbb{R}$. 
Changing the variables yields
$$
\lim_{\varepsilon\searrow 0}
\int_{-\infty}^\infty 
\exp(-\varepsilon \lvert\xi\rvert^k+ix\xi+it\xi^k)d\xi
=
\begin{cases}
t^{-1/k} A_k(t^{-1/k}x) 
& \text{if\,\,} t> 0,
\\
(-t)^{-1/k} \overline{A_k(-(-t)^{-1/k}x)} 
& \text{if\,\,} t< 0, 
\end{cases}
$$
and then, 
\begin{equation}
\left\lvert
\lim_{\varepsilon\searrow 0}
\int_{-\infty}^\infty 
\exp(-\varepsilon \lvert\xi\rvert^k+ix\xi+it\xi^k)d\xi
\right\rvert
=
\lvert t\rvert^{-1/k} 
\lvert A_k(\lvert t\rvert^{-1/k}x\,\operatorname{sgn} t)
\rvert
\label{Ak_eq}
\end{equation}
for $t\in\mathbb{R}\setminus\{0\}$. 
Here, 
$\operatorname{sgn}$ denotes the signature function: 
$\operatorname{sgn} s=s/\lvert s\rvert$ 
if $s\in\mathbb{R}\setminus\{0\}$, 
$\operatorname{sgn} s=0$ if $s=0$. 
\par
The bounds of $A_k$ are the following. 
\begin{prop}
$$
\left\lvert
\lim_{\varepsilon\searrow 0}
\int_{-\infty}^\infty 
\exp(-\varepsilon\lvert\xi\rvert^k+ix\xi+it\xi^k)d\xi
\right\rvert
\leqslant 
12 \lvert t\rvert^{-1/k}
$$
for any $k_j\in\{2,3,4,\ldots\}$ and $t\in\mathbb{R}\setminus\{0\}$. 
\label{prop:Ak_bdd_2}
\end{prop}
In view of \eqref{Ak_eq}, 
Proposition~\ref{prop:Ak_bdd_2} immediately follows 
from the following. 
\begin{lem}
$\lvert A_k(x)\rvert \leqslant 12$ 
for any $k_j\in\{2,3,4,\ldots\}$. 
\label{lem:Ak_bdd}
\end{lem}
\begin{proof}
$\lvert A_2(x)\rvert=2\sqrt{\pi}$ 
follows from the well-known formula
$$
\lim_{\varepsilon\searrow 0}
\int_{-\infty}^\infty 
\exp(-\varepsilon \xi^2+ix\xi+it\xi^2)d\xi
=
\frac{2\sqrt{\pi}}{\sqrt{it}}e^{-ix^2/4t}.
$$
Moreover, 
$$
A_3(x)
=\frac{2\pi}{\sqrt[3]{3}}
\operatorname{Ai}(x/\sqrt[3]{3}), 
$$
where $\operatorname{Ai}$ denotes the Airy function. 
The boundedness of $\operatorname{Ai}$ yields that of $A_3$. 
\par
Now, 
we will observe the boundedness of the functions $A_k$. 
The proof depends on whether $k$ is even or odd. 
Let $\varepsilon>0$. 
First, 
suppose that $k$ is even. 
Set
$$
x^\star
=
\left(\frac{\lvert 1-x\rvert}{k}\right)^{1/(k-1)}
\operatorname{sgn}(1-x),
\,\, 
x_\star
=
\left(\frac{\lvert -1-x\rvert}{k}\right)^{1/(k-1)}
\operatorname{sgn}(-1-x). 
$$
Note that these numbers satisfy
$$
x+k(x^\star)^{k-1}=1,
\,\, 
x+k(x_\star)^{k-1}=-1, 
$$
and
\begin{equation}
0
\leqslant 
x^\star-x_\star
\leqslant 
\frac{2}{k^{1/(k-1)}}
\leqslant 
2. 
\label{eq:star_even}
\end{equation}
Integrating by parts on the intervals 
$(-\infty,x_\star]$ and $[x^\star,\infty)$, 
we have
\begin{align*}
{}&
  \int_{-\infty}^\infty 
  \exp(-\varepsilon s^k+ixs+is^k)ds
\\
&=
  \int_{x_\star}^{x^\star} 
  \exp(-\varepsilon s^k+ixs+is^k)ds
\\*
&  \quad
  {}-i\frac{\exp(-\varepsilon (x_\star)^k+ixx_\star+i(x_\star)^k)}
  {i\varepsilon k(x_\star)^{k-1}-1}
  +i\frac{\exp(-\varepsilon (x^\star)^k+ixx^\star+i(x^\star)^k)}
  {i\varepsilon k(x^\star)^{k-1}+1}
\\*
&\quad
  {}-ik(k-1)(1+i\varepsilon)
  \int_{-\infty}^{x_\star}
  \frac{s^{k-2}\exp(-\varepsilon s^k+ixs+is^k)}
  {(i\varepsilon ks^{k-1}+x+ks^{k-1})^2}ds 
\\*
&\quad
  {}-ik(k-1)(1+i\varepsilon)
  \int_{x^\star}^\infty
  \frac{s^{k-2}\exp(-\varepsilon s^k+ixs+is^k)}
  {(i\varepsilon ks^{k-1}+x+ks^{k-1})^2}ds. 
\end{align*}
Changing the variables $x+ks^{k-1}=\tilde s$ 
in the second and the third integrals, 
and using \eqref{eq:star_even}, 
we have
\begin{align*}
\left\lvert
  \int_{-\infty}^\infty 
  \exp(-\varepsilon s^k+ixs+is^k)ds
  \right\rvert
&\leqslant 
  2+\int_{x_\star}^{x^\star} ds
  +k(k-1)(1+\varepsilon)
  \int_{-\infty}^{x_\star}
  \frac{s^{k-2}}{(x+ks^{k-1})^2}ds
\\*
&\quad
  +k(k-1)(1+\varepsilon)
  \int_{x^\star}^\infty
  \frac{s^{k-2}}{(x+ks^{k-1})^2}ds
\\
&= 
  2+x^\star-x_\star
  +(1+\varepsilon)
  \int_{-\infty}^{-1}
  \frac{d\tilde s}{\tilde s^2}
  +(1+\varepsilon)
  \int_1^\infty
  \frac{d\tilde s}{\tilde s^2}
\\
&= 
  2+x^\star-x_\star+2(1+\varepsilon)
\\
&\leqslant 
  6+2\varepsilon. 
\end{align*}
Therefore, 
we obtain $\lvert A_k(x)\rvert\leqslant 6$. 
This completes the proof in the case where $k$ is even. 
\par
Second, 
suppose that $k$ is odd. 
Set
$$
x^\star
=
\left(\frac{1-x}{k}\right)^{1/(k-1)}
\text{\,\,for\,\,}x\leqslant 1,
\,\, 
x_\star
=
\left(\frac{-1-x}{k}\right)^{1/(k-1)}
\text{\,\,for\,\,}x\leqslant -1
$$
this time. 
Note that these numbers satisfy
$$
x+k(\pm x^\star)^{k-1}=1,
\,\, 
x+k(\pm x_\star)^{k-1}=-1, 
$$
and
$$
x^\star
\leqslant 
\left(\frac{2}{k}\right)^{1/(k-1)}
\leqslant 
1,
\,\, 
0
\leqslant 
x^\star-x_\star
\leqslant 
\left(\frac{2}{k}\right)^{1/(k-1)}
\leqslant 
1. 
$$
When $x\leqslant -1$, 
integrating by parts on the intervals $(-\infty,-x^\star]$, 
$[-x_\star,x_\star]$ and $[x^\star,\infty)$ yields 
$\lvert A_k(x)\rvert\leqslant 12$. 
When $-1\leqslant x\leqslant 0$, 
integrating by parts on the intervals 
$(-\infty,-x^\star]$ and $[x^\star,\infty)$ yields 
$\lvert A_k(x)\rvert\leqslant 6$. 
When $0\leqslant x\leqslant 1$, 
integrating by parts on the intervals 
$(-\infty,-1]$ and $[1,\infty)$ yields 
$\lvert A_k(x)\rvert\leqslant 10/3$. 
When $x\geqslant 1$, 
integrating by parts on the whole interval yields 
$\lvert A_k(x)\rvert\leqslant 2$. 
This completes the proof in the case where $k$ is odd. 
\end{proof}
\par
Incidentally,
we state an estimate of another oscillatory integral.
\begin{cor}
$$
\left\lvert
\lim_{\varepsilon\searrow 0}
\int_{-\infty}^\infty 
\exp(-\varepsilon\lvert\xi\rvert^K+ix\xi+it\lvert\xi\rvert^K)d\xi
\right\rvert
\leqslant 
10\lvert t\rvert^{-1/K}. 
$$
for any $K>1$ and $t\in\mathbb{R}\setminus\{0\}$. 
\label{cor:Ak_bdd_noninteger}
\end{cor}
\par
We can verify Corollary~\ref{cor:Ak_bdd_noninteger} 
by similar calculus to the proof of Lemma~\ref{lem:Ak_bdd} 
in the case where $k$ is even. 
\par
Now, 
we prove the Fourier restriction theorem. 
\begin{proof}[Proof of Theorem~\ref{thm:main_thm}]
We argue as in \cite[Proof of Theorem~1, Case~I]{Strichartz1}. 
Set 
$$
R_0(\xi')
=
\sum_{j=1}^{n-1}
\lvert \xi_j\rvert^{k_j},
\,\,
\frac{1}{q}
=
\sum_{j=1}^{n-1}
\frac{1}{k_j}
$$
for short. 
Using Lemma~\ref{lem:formula_01}, 
we have
\begin{align*}
{}&
  \check G_z(x)
\\
&=
  (2\pi)^{-n/2}
  \lim_{\varepsilon\searrow 0}
  \int_{\mathbb{R}^n}
  e^{-\varepsilon (R_0(\xi')+\lvert \xi_n-R(\xi')\rvert)}
  G_z(\xi)e^{ix\cdot\xi}d\xi
\\
&=
  (2\pi)^{-n/2}
  \lim_{\varepsilon\searrow 0}
  \int_{\mathbb{R}^n}
  \frac{e^{-\varepsilon (R_0(\xi')+\lvert \xi_n-R(\xi')\rvert)}
  (\xi_n-R(\xi'))_+^z e^{ix\cdot\xi}}
  {\Gamma(z+1)}d\xi
\\
&=
  (2\pi)^{-n/2}
  \lim_{\varepsilon\searrow 0}
  \int_{\mathbb{R}^{n-1}}
  \exp(-\varepsilon R_0(\xi')+i x'\cdot\xi')
  \int_{-\infty}^\infty 
  \frac{e^{-\varepsilon \lvert \xi_n-R(\xi')\rvert}
  (\xi_n-R(\xi'))_+^z e^{ix_n\xi_n}}
  {\Gamma(z+1)}d\xi_n d\xi'
\\
&=
  (2\pi)^{-n/2}
  \lim_{\varepsilon\searrow 0}
  \int_{\mathbb{R}^{n-1}}
  \exp(-\varepsilon R_0(\xi')+ix_nR(\xi')+i x'\cdot\xi')
  \int_0^\infty 
  \frac{e^{-\varepsilon \xi_n}\xi_n^z e^{ix_n\xi_n}}
  {\Gamma(z+1)}d\xi_nd\xi'
\\
&=
  i(2\pi)^{-n/2}
  e^{iz\pi/2}(x_n+i0)^{-z-1}
  \lim_{\varepsilon\searrow 0}
  \int_{\mathbb{R}^{n-1}}
  \exp{(-\varepsilon R_0(\xi')+ix_nR(\xi')
  +i x'\cdot\xi')}d\xi'
\\
&=
  i(2\pi)^{-n/2}
  e^{iz\pi/2}(x_n+i0)^{-z-1}
  \prod_{j=1}^{n-1}
  \lim_{\varepsilon\searrow 0}
  \int_{-\infty}^\infty 
  \exp(-\varepsilon\lvert \xi_j\rvert^{k_j}
  +ix_j\xi_j+ia_jx_n \xi_j^{k_j})d\xi_j. 
\end{align*}
For $z\in\mathbb{C}$, 
$\Re z$ and $\Im z$ denote its 
real part and imaginary part, 
respectively.
Therefore, 
using Proposition~\ref{prop:Ak_bdd_2}, 
we have
\begin{align*}
  \lvert\check G_z(x)\rvert
&=
  (2\pi)^{-n/2}
  e^{-\Im z\pi/2}
  \lvert x_n\rvert^{-\Re z-1}
  \prod_{j=1}^{n-1}
  \left\lvert
  \lim_{\varepsilon\searrow 0}
  \int_{-\infty}^\infty 
  \exp(-\varepsilon\lvert\xi_j\rvert^{k_j}
  +ix_j\xi_j+ia_jx_n \xi_j^{k_j})d\xi_j
  \right\rvert
\\
&\leqslant 
  12^{n-1}(2\pi)^{-n/2}
  e^{-\Im z\pi/2}
  \lvert x_n\rvert^{-\Re z-1-1/q}
  \prod_{j=1}^{n-1}
  \lvert a_j\rvert^{1/k_j}
\end{align*}
for $x_1,\ldots,x_n\neq 0$. 
Namely, 
we obtain
$$
\lvert
\check G_{-(1+1/q)+i\eta}(x)
\rvert
\leqslant 
Ce^{-\eta\pi/2}
$$
for $x_1,\ldots,x_n\neq 0$ and all $\eta\in\mathbb{R}$ 
with $C$ depending only on $n$ and $a_1,\ldots,a_{n-1}$. 
Now, 
we can apply Proposition~\ref{prop:Strichartz_method} 
with $\lambda=1+1/q$. 
Then we obtain the desired sufficient condition 
on the exponent $p$ for \eqref{eq:restriction_ineq_rewritten}. 

In the rest of the proof, 
we also argue in essentially the same way as 
\cite[Proof of Theorem~1]{Strichartz1}. 
We use a homogeneity argument 
with respect to the nonisotropic dilations
$$
d_s \phi(x)
=\phi(s^{1/k_1}x_1,\ldots ,s^{1/k_{n-1}}x_{n-1},sx_n)
$$
for $s>0$. 
On one hand, 
\begin{align*}
  (d_s \phi)^\wedge (\xi)
&=
  (2\pi)^{-n/2}
  \int_{\mathbb{R}^n}
  \phi(s^{1/k_1}x_1,\ldots ,s^{1/k_{n-1}}x_{n-1},sx_n)
  e^{-ix\cdot \xi}dx
\\
&=
  (2\pi)^{-n/2}
  s^{-(1+1/q)}
  \int_{\mathbb{R}^n}\phi(x)
\\*
&\quad{}\times
  \exp(-is^{-1/k_1}x_1\xi_1-\cdots-is^{-1/k_{n-1}}x_{n-1}\xi_{n-1}
  -is^{-1}x_n\xi_n)dx
\\
&=
  s^{-(1+1/q)}
  d_{s^{-1}}\hat \phi(\xi). 
\end{align*}
On the other hand, 
$$
\int_{\mathbb{R}^{n-1}} 
\lvert
d_{s^{-1}}\hat \phi(\xi',R(\xi'))
\rvert^2d\xi'
=
s^{1/q}
\int_{\mathbb{R}^{n-1}} 
\lvert
\hat \phi(\xi',R(\xi'))
\rvert^2
d\xi'. 
$$
Now using
$$
\lVert
d_s \phi
\rVert_{L^p(\mathbb{R}^n)}
=
s^{-(1+1/q)/p}
\lVert\phi\rVert_{L^p(\mathbb{R}^n)}, 
$$
and applying \eqref{eq:restriction_ineq_rewritten} 
for the function $d_s \phi$, 
we have
\begin{align*}
  \left(\int_{\mathbb{R}^{n-1}} 
  \lvert
  \hat \phi(\xi',R(\xi'))
  \rvert^2
  d\xi'\right)^{1/2}
&= 
  s^{-1/2q}
  \left(
  \int_{\mathbb{R}^{n-1}} 
  \lvert
  d_{s^{-1}}\hat \phi(\xi',R(\xi'))
  \rvert^2
  d\xi'
  \right)^{1/2}
\\
&= 
  s^{1+1/2q}
  \left(
  \int_{\mathbb{R}^{n-1}} 
  \lvert
  (d_s \phi)^\wedge(\xi',R(\xi'))
  \rvert^2
  d\xi'
  \right)^{1/2}
\\
&\leqslant 
  C_p
  s^{1+1/2q}
  \lVert
  d_s \phi
  \rVert_{L^p(\mathbb{R}^n)}
\\
&= 
  C_p
  s^{1+1/2q-(1+1/q)/p}
  \lVert\phi\rVert_{L^p(\mathbb{R}^n)}. 
\end{align*}
Therefore, 
to obtain \eqref{eq:restriction_ineq_rewritten} 
for any $s>0$, 
we must have 
$1+1/2q-(1+1/q)/p=0$. 
Then the desired necessary condition 
on the exponent $p$ for 
\eqref{eq:restriction_ineq_rewritten} follows. 
This completes the proof. 
\end{proof}
\par
In view of Corollary~\ref{cor:Ak_bdd_noninteger}, 
we can prove Theorem~\ref{thm:main_thm} 
with replacing some $\xi_j^{k_j}$ 
in \eqref{eq:Rxi_prime} by 
$\lvert\xi_j\rvert^{K_j}$ where $K_j>1$. 
\par
For the case where $R$ is homogeneous 
($k_1=\cdots=k_{n-1}$), 
see \cite[Chapter~8, \S~5.17]{temp003}. 
\section{An application to partial differential equations}
\label{sect:pde}
Finally, 
we prove Theorem~\ref{thm:main_cor}. 
\par
\begin{proof}[Proof of Theorem~\ref{thm:main_cor}]
We argue as in \cite[Proof of Corollary~1]{Strichartz1}. 
We may assume 
$\phi\in\mathscr{S}(\mathbb{R}^n)$ 
and $f\in\mathscr{S}(\mathbb{R}^{1+n})$, 
where $\mathscr{S}$ denotes the Schwartz class.
Set 
$$
R_0(\xi)
=
\sum_{j=1}^n 
\lvert\xi_j\rvert^{k_j},
\,\,
\frac{1}{q}
=
\sum_{j=1}^n
\frac{1}{k_j}
$$
for short. 
The solution to 
\eqref{eq:pde_1}--\eqref{eq:pde_2} is written as
\begin{align*}
  u(t,x)
&=
  (2\pi)^{-n/2}
  \int_{\mathbb{R}^n}
  e^{ix\cdot \xi}
  e^{ita(\xi)}
  \hat \phi(\xi)
  d\xi
\\*
&\quad{}+
  (2\pi)^{-n/2}i
  \int_0^t 
  \lim_{\varepsilon\searrow 0}
  \int_{\mathbb{R}^n}
  e^{-\varepsilon R_0(\xi)}
  e^{ix\cdot \xi}
  e^{i(t-s)a(\xi)}
  (f(s,\cdot))^\wedge(\xi)
  d\xi ds. 
\end{align*}
By duality, 
\eqref{eq:restriction_ineq_rewritten} is equivalent to
\begin{align*}
  \lVert
  (F d\mu_n)^\wedge
  \rVert_{L^{p/(p-1)}(\mathbb{R}^n)}
\leqslant 
  C_{p/(p-1)}
  \left(
  \int_{\mathbb{R}^{n-1}}
  \lvert 
  F(x',\tilde R(\xi'))
  \rvert^2 
  d\xi'
  \right)^{1/2}
  \text{\,\,for\,\,}F \in L^2(d\mu_n)&. 
\end{align*}
Now, 
let $S=\{(t,\xi)\in\mathbb{R}^{1+n};\,t+a(\xi)=0\}$. 
Replacing $n$ by $n+1$ and $x_{n+1}$ by $t$, 
we have
$$
  \lVert
  \mathscr{F}_{t,x}[F d\mu_t]
  \rVert_{L^{p/(p-1)}(\mathbb{R}^{1+n})}
\leqslant 
  C_{p/(p-1)}
  \left(
  \int_{\mathbb{R}^n}
  \lvert 
  F(-a(x),x)
  \rvert^2 
  dx
  \right)^{1/2}
  \text{\,\,for\,\,}F \in L^2(d\mu_t), 
$$
where $\mathscr{F}_{t,x}[f]$ denotes 
the Fourier transform of $f$ in $(t,x)\in\mathbb{R}^{1+n}$,
that is,
$$
\mathscr{F}_{t,x}
[f](\tau,\xi)
=
(2\pi)^{-(1+n)/2}
\int_{\mathbb{R}^n}
\int_{-\infty}^\infty 
f(t,x)
e^{-it\tau-ix\cdot\xi}
dtdx. 
$$
Applying this, 
we have
\begin{align*}
  \left\lVert
  \int_{\mathbb{R}^n}
  e^{ix\cdot \xi}
  e^{ita(\xi)}
  \hat \phi(\xi)
  d\xi
  \right\rVert_{L^{p/(p-1)}(\mathbb{R}^{1+n})}
&=
  \left\lVert
  \int_S 
  e^{ix\cdot \xi}
  e^{-it\tau}
  \hat \phi(\xi)
  d\mu_t(\tau,\xi)
  \right\rVert_{L^{p/(p-1)}(\mathbb{R}^{1+n})}
\\
&=
  (2\pi)^{(1+n)/2}
  \lVert
  \mathscr{F}_{\tau,\xi}[\hat \phi(\xi)d\mu_t](t,-x)
  \rVert_{L^{p/(p-1)}(\mathbb{R}^{1+n})}
\\
&\leqslant 
  C
  \left(
  \int_{\mathbb{R}^n}
  \lvert \hat \phi(x)\rvert^2
  dx
  \right)^{1/2}
\\
&=
  C\lVert\phi\rVert_{L^2(\mathbb{R}^n)}. 
\end{align*}
Next, 
set
$$
T_sf(x)
=
\lim_{\varepsilon\searrow 0}
\int_{\mathbb{R}^n}
e^{-\varepsilon R_0(\xi)}
e^{ix\cdot \xi}
e^{i(t-s)a(\xi)}
\hat f(s,\xi)d\xi
$$
for $s\in\mathbb{R}$. 
Here, 
we write 
$\hat f(s,\xi)=(f(s,\cdot))^\wedge(\xi)$.
It follows from the Plancherel theorem that
\begin{equation}
\lVert T_sf\rVert_{L^2(\mathbb{R}^n)}
=
\lVert f(s,\cdot)\rVert_{L^2(\mathbb{R}^n)}. 
\label{eq:2_2bdd}
\end{equation}
Now, 
we have
\begin{align*}
{}&
  T_sf(x)
\\
&=
  (2\pi)^{-n/2}
  \lim_{\varepsilon\searrow 0}
  \int_{\mathbb{R}^n}
  e^{-\varepsilon R_0(\xi)}
  e^{ix\cdot \xi}
  e^{i(t-s)a(\xi)}
  \int_{\mathbb{R}^n}
  f(s,y)
  e^{-iy\cdot \xi}
  dyd\xi
\\
&=
  (2\pi)^{-n/2}
  \int_{\mathbb{R}^n}
  f(s,y)
  \left(
  \prod_{j=1}^n 
  \lim_{\varepsilon\searrow 0}
  \int_{-\infty}^\infty 
  \exp(-\varepsilon\lvert\xi\rvert^{k_j}
  +i(x_j-y_j)\xi_j+ia_j(t-s)\xi_j^{k_j})d\xi_j
  \right)
  dy. 
\end{align*}
Using Proposition~\ref{prop:Ak_bdd_2}, 
we have
\begin{align*}
  \lvert T_sf(x)\rvert
&\leqslant 
  (2\pi)^{-n/2}
  \int_{\mathbb{R}^n}
  \lvert f(s,y)\rvert
  \prod_{j=1}^n 
  \left\lvert
  \int_{-\infty}^\infty 
  \exp(i(x_j-y_j)\xi_j+ia_j(t-s)\xi_j^{k_j})
  d\xi_j
  \right\rvert dy
\\
&\leqslant 
  12^n (2\pi)^{-n/2} 
  \lvert t-s\rvert^{-1/q}
  \lVert f(s,\cdot)\rVert_{L^1(\mathbb{R}^n)}
  \prod_{j=1}^n 
  \lvert a_j\rvert^{-1/k_j}
\end{align*}
for $s\neq t$. 
Namely, 
\begin{equation}
\lVert T_sf\rVert_{L^\infty(\mathbb{R}^n)}
\leqslant
C\lvert t-s\rvert^{-1/q}
\lVert f(s,\cdot)\rVert_{L^1(\mathbb{R}^n)}. 
\label{eq:infty_1bdd}
\end{equation}
Interpolating \eqref{eq:2_2bdd}--\eqref{eq:infty_1bdd}, 
we have
\begin{align*}
  \lVert T_sf\rVert_{L^{p/(p-1)}(\mathbb{R}^n)}
&\leqslant 
  C
  \lvert t-s\rvert^{-(2/p-1)/q}
  \lVert f(s,\cdot)\rVert_{L^p(\mathbb{R}^n)}\\
&=
  C\lvert t-s \rvert^{-2(p-1)/p}
  \lVert f(s,\cdot)\rVert_{L^p(\mathbb{R}^n)}. 
\end{align*}
Now, 
using the Hardy-Littlewood-Sobolev inequality 
(see \cite[page 354]{temp003} for instance), 
we have
\begin{align*}
{}&\left\lVert\int_0^t
  \int_{\mathbb{R}^n}
  e^{ix\cdot \xi}
  e^{i(t-s)a(\xi)}
  \hat f(s,\xi)
  d\xi ds
  \right\rVert_{L^{p/(p-1)}(\mathbb{R}^{1+n})}
\\
&=
  \left\lVert
  \int_0^t T_sf(x) ds
  \right\rVert_{L^{p/(p-1)}(\mathbb{R}^{1+n})}
\\
&\leqslant
  \left\lVert
  \int_0^t 
  \lVert T_sf\rVert_{L^{p/(p-1)}(\mathbb{R}^n)}
  ds
  \right\rVert_{L^{p/(p-1)}(\mathbb{R}_t)}
\\
&\leqslant
  \left\lVert
  \int_{-\infty}^\infty 
  \lVert T_sf\rVert_{L^{p/(p-1)}(\mathbb{R}^n)} 
  ds
  \right\rVert_{L^{p/(p-1)}(\mathbb{R}_t)}
\\
&\leqslant
  C
  \left\lVert
  \int_{-\infty}^\infty 
  \lvert t-s\rvert^{-2(p-1)/p}
  \lvert f(s,\cdot)\rVert_{L^p(\mathbb{R}^n)}
  ds
  \right\rVert_{L^{p/(p-1)}(\mathbb{R}_t)}
\\
&\leqslant 
  C
  \lVert\,\lVert f(t,\cdot)
  \rVert_{L^p(\mathbb{R}^n)}\rVert_{L^p(\mathbb{R}_t)}
\\
&=
  C 
  \lVert f\rVert_{L^p(\mathbb{R}^{1+n})}. 
\end{align*}
Thus, 
we obtain
$$
\lVert u\rVert_{L^{p/(p-1)}(\mathbb{R}^{1+n})}
\leqslant
C
(\lVert\phi\rVert_{L^2(\mathbb{R}^n)}
+\lVert f\rVert_{L^p(\mathbb{R}^{1+n})}). 
$$
This completes the proof. 
\end{proof}
\section*{Acknowledgment}
The author expresses gratitude to Hiroyuki Chihara 
for helpful advices. 

\bibliographystyle{amsplain}

\end{document}